%% file: main.tex
\documentclass[journal]{IEEEtran} 

\usepackage{cite}
\usepackage[]{graphicx}
\graphicspath{{figures/}}

\usepackage{amsmath}
  {

	\interdisplaylinepenalty=2500
  }
\usepackage{algorithm}
\usepackage[noend]{algpseudocode}
\usepackage{amsfonts}
\usepackage{amssymb}
\usepackage{booktabs}
\usepackage[caption=false,font=footnotesize]{subfig}
\usepackage{stfloats}
\usepackage{enumitem}
\usepackage{tikz}
\usetikzlibrary{shapes,arrows}

\usepackage{hyperref}[6.83]
\hypersetup{colorlinks,
			linkcolor=[rgb]{.61,0,0.3},
			citecolor=[rgb]{.14,.47,.14}}

\usepackage{tikz}
\usetikzlibrary{arrows.meta,
	calc, chains,
	quotes,
	positioning,
	shapes.geometric}

\newcommand{\norm}[1]{\left\Vert #1 \right\Vert}

\definecolor{OliveGreen}{rgb}{0.33, 0.42, 0.18}

\makeatletter
\newcommand{\removelatexerror}{\let\@latex@error\@gobble}
\makeatother

\newcommand{\maxx}{\mathop{\bf max}}

\definecolor{green2}{rgb}{0.0, 0.5, 0.0}

\newcommand{\review}[1]{{#1}}

\begin{document}

\title{Efficient Computation of Power System Maximum Transient Linear Growth}
\author{Daniel Adrian Maldonado,~\IEEEmembership{Member,~IEEE,}
        ~Emil~M.~Constantinescu,~\IEEEmembership{Member,~IEEE}
        Junbo Zhao,~\IEEEmembership{Senior Member,~IEEE,}
        and~Mihai Anitescu,~\IEEEmembership{Member,~IEEE}
\thanks{This material is based upon work supported by the U.S. Department of Energy Advanced Grid Modernization program.}
\thanks{Authors are with the Mathematics and Computer Science Division at 
Argonne National Laboratory, Lemont, Illinois, U.S.A.,and with the University of Connecticut, Storrs, Connecticut, U.S.A.}}%

\markboth{}%
{Analysis of small perturbation growth in power system dynamics}

\maketitle
\begin{abstract}
Existing methods to determine the stability of a power system to small perturbations are based on eigenvalue analysis and focus on the asymptotic (long-term) behavior of the power grid.  During the preasymptotic (short-term) transient, however, the system can exhibit large growth that is not explained by eigenvalues alone. In this paper we propose  a new framework to determine the maximum (optimal) preasymptotic growth using the singular value decomposition. The approach is tailored to the analysis of quantities of interest in power system dynamics, such as the set of rotor speed deviations. Matrix-free techniques are developed to avoid the explicit formation of dense matrices and enable the analysis of large-scale systems without reaching memory bounds. Extensive results carried out from small to very large-scale systems (e.g., 70k-bus system) verify the theoretical aspects of the technique.
\end{abstract}


\IEEEpeerreviewmaketitle 

\section{Introduction}
\label{section:intro}
\input{sections/intro} 

\section{Limitations of Small-Signal Stability Analysis}
\label{section:modal}
\input{sections/modal} 

\section{Proposed Method}
\label{section:methodology}
\input{sections/methodology} 

\section{Numerical Results}
\label{section:casestudies}
\input{sections/cases}

\section{Conclusions} \label{section:conclusion}
In this paper we have shown that, while small signal analysis based on eigenvalues can provide an understanding of the asymptotic behavior of the power system dynamics around an equilibrium, the preasymptotic response is strongly dependent on the condition of the eigenbasis. We have introduced a method to examine this behavior using  SVD, and we have tailored it to the analysis of power systems. We have introduced a matrix-free approach that allows us to analyze systems of large scale without being bounded by the available memory. Our test cases  illustrate how non-normality influences the preasymptotic response. Future work will extend this approach beyond the study of equilibrium perturbations and will incorporate the response of the system to stochastic forcing.

\section*{Acknowledgment}
This material was based upon work supported by the U.S. Department of Energy Advanced Grid Modernization program.
\bibliographystyle{ieeetr}
\bibliography{bibliography}

\vfill
\begin{flushright}
{\footnotesize
  \framebox{\parbox{0.5\textwidth}{
The submitted manuscript has been created by UChicago Argonne, LLC, 
Operator of Argonne National Laboratory (``Argonne"). Argonne, a
U.S. Department of Energy Office of Science laboratory, is operated
under Contract No. DE-AC02-06CH11357. The U.S. Government retains for
itself, and others acting on its behalf, a paid-up nonexclusive,
irrevocable worldwide license in said article to reproduce, prepare
derivative works, distribute copies to the public, and perform
publicly and display publicly, by or on behalf of the Government.
The Department of
Energy will provide public access to these results of federally sponsored research in accordance
with the DOE Public Access Plan. http://energy.gov/downloads/doe-public-access-plan. }}
\normalsize
}
\end{flushright}

\end{document}

%% file: sections/intro.tex
Owing to the increase of distributed energy resources such as photovoltaic  generators, small-scale battery energy storage systems, and other inverter-based resources, the presence of harmonics and other types of perturbations in the grid is increasing. In order to accommodate these resources,  a proper analysis is essential to determine their impact on the operation of the grid. Of particular importance is whether some of these small perturbations have the potential to grow into oscillations or other transients that could harm the grid.

These types of studies have been commonplace in power systems research, and as a result  several methods have been developed to study the effect of small perturbations in power systems \cite{Milano2020}. Some of these consist of the well-known modal analysis techniques that were developed between the 1960s and 1970s; and today, thanks to advances that allow efficient computation of the eigenvalues of large sparse matrices, these techniques are employed by a wealth of commercial transient stability  packages \cite{undrill, Anderson1969, Wong1988}.

Modal small-signal stability analysis techniques continue  to be used for myriad applications. For instance,  modal analysis is used to optimize the dynamic performance of low-inertia systems \cite{Venkatraman2021}. Modal constraints have also been embedded in the optimal power flow to account for small-signal stability \cite{Liu2022}. Machine learning and modal analysis have been combined to detect events from synchrophasor data \cite{Taghipourbazargani2022}. A higher-order extension to modal analysis based on the normal form also has been investigated \cite{Ugwuanyi2020}. These recent examples show  that research using modal analysis techniques is far from being exhausted.

Modal analysis, however, is not without its problems. Because of linearization, one of the main issues is that it is difficult to quantify what a small perturbation \textit{is} without relying on more advanced global techniques \cite{Alberto2012}. Another limitation of modal analysis is that eigenvalues characterize the asymptotic behavior of the underlying system.  During the initial period, however, a stable linear system can exhibit transient growth that is not desirable from an engineering point of view.

A shift toward alternatives to modal analysis began in the context of analyzing the transition to turbulence in fluid dynamics. In the 1990s, Baggett, Driscoll, and Trefethen observed that in a stable linearized model of fluid flow, a transition to turbulence could be seen due to large linear growth exciting the nonlinear part of the equations \cite{Baggett1995}. This large linear growth was attributed to the non-normality of the Jacobian matrix. Later, Farrell and Ioannou devised a method to characterize such non-normal growth by using the singular value decomposition (SVD) of the exponential map \cite{Farrell1996}. Another way of studying this non-normal growth was introduced by Trefethen and Embree by the use of pseudospectral analysis \cite{Trefethen2005}. These techniques have been termed non-modal, in contrast to traditional modal analysis that uses the eigenvalues of the linearized system.

The focus of this paper is the development of a non-modal technique based on  SVD to analyze small perturbations in power system dynamics and determine the largest possible transient growth. We  show how this technique can offer additional insight into  the preasymptotic dynamic behavior of power systems. The main contributions are summarized as follows.

\begin{itemize}
    \item Non-modal techniques are developed to analyze power system dynamics. We develop specific techniques to characterize the preasymptotic growth of a subset of variables, making it applicable to traditional transient stability analysis, such as rotor stability.
    \item A matrix-free approach is developed to compute the SVD of the exponential operator for the differential-algebraic equations (DAEs). Our approach avoids the explicit formation of dense matrices, thus making the technique scalable to large-scale power systems.
\end{itemize}

Our research connects with recent evolving efforts in the physics community to study network structure, its effect on stability, and the concept of \textit{reactivity} \cite{Asllani2018, Muolo2020, Duan2022}. Physics researchers have shown that network reactivity, defined as the capacity of the linear system to amplify its responses to small perturbations, is linked to the structural properties of the network and the presence of non-normality \cite{Duan2022}. This non-normality can result in high sensitivity of the eigenvalues to noise and can reduce their physical meaning \cite{Asllani2018}.

Our approach shares some ideas with $H_{\infty}$ and $H_{2}$ robust control methods in power systems  \cite{Chaudhuri2004, Zhao2018}. $H_{\infty}$ and $H_{2}$ methods make extensive use of matrix norms and the SVD. However, these methods compute norms in the Hardy space and usually focus on analyzing the transfer function in the frequency domain. While there are time domain interpretation of such norms, they are always defined in terms of integrals in time, and as a result do not quantify well the possibly large system response over a small time interval. In contrast, we are concerned with the matrix norm of the exponential operator induced by a vector norm, and its behavior over a finite interval of time.

The paper begins with a short discussion on the tenets behind the small-signal analysis, and we show how conditioning and non-normality play an important role in the preasymptotic behavior of the linearized system. In Section \ref{section:methodology} we introduce  SVD to compute the maximum preasymptotic growth of the system. We also discuss the importance of choosing appropriate norms, and we detail how to compute efficiently. Section \ref{section:casestudies} develops case studies to illustrate the basic premises of the technique and how it can be applied to detailed, large-scale systems. Section V summarizes our work and briefly presents future extensions. 

%% file: sections/modal.tex
To analyze the effect of small perturbations on the stability of a nonlinear dynamical system, Lyapunov proposed a technique that is now called Lyapunov's first method or Lyapunov's indirect method. Given a nonlinear dynamical system $\dot{x} = f(x)$, its linearization around a steady state $x_o$ for a small perturbation $\delta x(t)$ can be written as
\begin{equation}
    \delta \dot{x}(t) = f(x_o + \delta x(t)) = A\delta x(t) + O(||\delta x(t)||^2) \,.
\end{equation}
Here we consider $x \in \mathbf{R}^N$, $x_o$ is such that $f(x_o) = 0$, and $A$ is the Jacobian matrix of $f()$. If the perturbation $\delta x(t)$ is small enough that the higher-order term is negligible, we can say that the motion of the system is governed by the linear operator $A$. In this case the system will be stable if the eigenvalues of $A$ fall on the left side of the complex plane. The presence of positive eigenvalues will result, if excited, in exponential growth or undamped oscillatory behavior. Denoting $u(t) = \delta x(t)$, the linear dynamical system that governs the small-signal response of the system is
\begin{equation}
    \dot{u}(t) = Au(t) \,,
\end{equation}
and the solution:
\begin{equation}
    u(t) = e^{At}u(0) \,.
\end{equation}
If $A$ is diagonalizable, the set of eigenvectors $v_i$ of $A$ forms a basis of $\mathbf{R}^N$, and we can write
\begin{equation}
    u(0) = \alpha_1 v_1 + \alpha_2 v_2 + \dots + \alpha_n v_2 \,,
\end{equation}
where $v_i$ are the eigenvectors of $A$ and $\alpha_i$ the scalar coefficients of the basis representation. Since $e^{At}$ is defined via the infinite power series
\begin{equation}
    e^{At} = I + At + \frac{1}{2!}A^2t^2 + \frac{1}{3!}A^3t^3 + \dots \,,
\end{equation}
and since $A^kv_i=\lambda^k_i v_i$ with $\lambda_i$ the eigenvalue associated to $v_i$,  we have that
\begin{align}
\label{eq:modalresp}
\begin{split}
   u(t) = e^{At}u(0) & = \sum_{k=0}^N \alpha_k (1 + \lambda_kt + \lambda^2_k t^2 + \dots)v_k \,\\
   & =  \sum_{k=0}^N \alpha_k e^{\lambda_k t} v_k \,.
\end{split}
\end{align}
In the limit when $t \rightarrow \infty$, the magnitude of $\norm{u(t)}$ is dominated by the eigenvalue with the largest real part, and it will decrease to zero if the real part of all the eigenvalues is negative. By examining the eigenvalues of the Jacobian matrix, small-signal stability analysis makes \textit{asymptotic} (i.e., as $ t \rightarrow \infty $) assertions about the behavior of the system under small disturbances of the steady state.

In this paper we are concerned with transient preasymptotic growth and its consequences on the stability of the power system. One can see that in the preasymptotic regime, the eigenvalues alone do not paint a complete picture of the system behavior. Given a linear system defined by the operator $A$, the norm of the state at time $t$ can be bounded \cite{Trefethen2005-vv}:
\begin{align}
\label{eq:condition}
\begin{split}
   \norm{u(t)} = \norm{e^{At} u(0)} &\leq \norm{e^{At}}\norm{u(0)} \\
   &\leq \norm{Ve^{t\Lambda} V^{-1}} \norm{u(0)} \\
   &\leq \norm{V}\norm{V^{-1}} \maxx_{\lambda \in \sigma (A)} \lvert e^{t\lambda} \rvert \norm{u(0)} \\
   &\leq \kappa (V) \maxx_{\lambda \in \sigma (A)} \lvert e^{t\lambda} \rvert \norm{u(0)} \,.
\end{split}
\end{align}
For simplicity, we assume $A$ is diagonalizable. Here $\Lambda$ is the diagonalized operator, $V$ is the change of basis matrix, $\sigma(A)$ is the spectrum of $A$, and $\kappa(V)$ is the condition number of $V$. Only in the case where $\kappa(V)$ is equal to $1$, which is true when $V$ is orthogonal, will the behavior of the linear system  be completely determined by its largest eigenvalue. One cannot  always find an orthogonal eigenbasis, except in the case where the operator $A$ is normal,
\begin{equation}
    A^*A = AA^* \,.
\end{equation}

Despite the asymptotic stability assurances of the modal approach, in the short term the transient growth is strongly dependent on the conditioning of the eigenbasis of the linear operator, and it can result in large transient growth when $\kappa(V)$ is large. In the following section we describe how to characterize this preasymptotic growth, and we develop an algorithm to compute the perturbation direction that induces it.

%% file: sections/methodology.tex
To examine the small-signal response of linearized power systems, we focus on the preasymptotic transient growth phenomenon. We can infer, from results such as the one shown in \eqref{eq:condition}, that in a system where this non-normality exists, small perturbations along specific directions can be amplified to such an extent that the system becomes unstable. From this perspective, analyzing the magnitude of this transient growth and the perturbation direction that produces it can provide additional insights into using traditional modal approaches.

\subsection{System Energy and Growth Function}
Farrell and Ioannou introduced the concept of system energy and growth function in the context of non-normal dynamics and subsequent analysis of growth \cite{Farrell1996}. Given a linear system
\begin{equation}
    \dot{x} = Ax \,, \quad x(0) = x_0 \,,
\end{equation}
we can define the energy of the system to be the function
\begin{equation}
    \mathcal{E}(t) = \norm{x(t)}^2 \,,
\end{equation}
and we define the growth as the energy at time $t$ relative to the initial conditions
\begin{equation}
    \mathcal{G}(t) = \frac{\mathcal{E}(t)}{\norm{x(0)}^2}  = \frac{\norm{x(t)}^2}{\norm{x(0)}^2}\,.
\end{equation}
The optimal growth is then defined as an optimization problem:
\begin{equation}
\label{eq:opt_problem}
\begin{aligned}
& \underset{x(0)}{\text{maximize}}
& & \mathcal{G}(t) \\
& \text{subject to}
& & \norm{x(0)} = 1 \,.
\end{aligned}
\end{equation}

If we use a Euclidean norm, we can see that
\begin{equation}
    \maxx_{\substack{\norm{x(0)}_2 \leq 1 \\ x(0) \neq 0}} \mathcal{G}(t) = \maxx_{\substack{\norm{x(0)}_2 \leq 1 \\ x(0) \neq 0}} \frac{ \norm{e^{At} x(0)}^2_2}{\norm{x(0)}^2_2} = \norm{e^{At}}^2_2 \,.
\end{equation}
Furthermore,  by computing the SVD of $e^{At}$, we get
\begin{equation}
    \norm{e^{At}}^2_2 = \norm{U\Sigma V^*}^2_2 = \norm{\Sigma}_2^2 = (\sigma_1(e^{At}))^2 \,,
\end{equation}
where $\sigma_1(A)$ is the largest singular value of $A$. We can then use the SVD to compute both the largest energy growth as well as the direction of the initial condition that causes it. This last direction can be found by the right singular vector associated with the largest singular value.

This framework is adequate in problems where the magnitudes of the different state variables are commensurate (e.g., fluid flow); but in electrical power systems a homogeneous measure can be misleading. To wit, the rotor speed and the subtransient voltage are vastly different in magnitude, and the latter can dominate the growth function. In addition, we are often concerned with the behavior of a group of quantities of interest, such as the collection of rotor speeds (rotor stability analysis) or voltage magnitudes (voltage stability analysis). For instance, given the linear system that arises from the linearization of the classical electromechanical model
\begin{equation}
\label{eq:classical_linear}
    A = \begin{bmatrix} 0 & I \\ -H^{-1}K & -H^{-1}D\end{bmatrix} \,, \quad
    x(t) = \begin{bmatrix} \delta  \\  \omega\end{bmatrix} \,,
\end{equation}
where $H$ and $D$ are diagonal matrices of machine inertia and damping, respectively, $K$ is the lossless admittance matrix, $\delta$ the rotor angle, and $\omega$ the rotor speeds, one often  wishes to define the energy as the sum of the kinetic and potential energy $ E_T = \frac{1}{2} \delta^\top  K \delta + \frac{1}{2} \omega^\top  H \omega$. Recall that, given a positive definite matrix $C \in \mathbf{C}^{n \times n}$, we can define an associated $C-$inner product with a vector norm $\norm{\cdot}_C$, and it can be shown that $\norm{x}_C = \norm{Cx}_2$. Then
\begin{equation}
    \maxx_{\substack{\norm{x(0)}_2 \leq 1 \\ x(0) \neq 0}} \frac{ \norm{e^{At} x(0)}^2_C}{\norm{x(0)}^2_2} = \maxx_{\substack{\norm{x(0)}_2 \leq 1 \\ x(0) \neq 0}} \frac{\norm{Ce^{At} x(0)}^2_2}{\norm{x(0)}^2_2} \,,
\end{equation}
where the last term is just the induced matrix norm $\norm{Ce^{At}}^2_2$.
If we use the energy norm, we would have 
\begin{equation}
\label{eq:energy_norm}
    C = \begin{bmatrix} K^{\frac{1}{2}} & \\ & H^{\frac{1}{2}} \end{bmatrix} \,,
\end{equation}
with $K^{\frac{1}{2}}$ and $H^{\frac{1}{2}}$ being the Cholesky factors of $K$ and $H$, respectively. More generally, we can define the subordinate matrix norm to be the vector norms $\norm{\cdot}_C$ and $\norm{\cdot}_B$ as 
\begin{equation}
    \norm{e^{At}}_{(C,B)} = \maxx_{\substack{\norm{x}_B \leq 1 \\ x \neq 0}} \frac{ \norm{e^{At} x}^2_C}{\norm{x}^2_B} = \maxx_{\substack{\norm{y}_2 \leq 1 \\ y \neq 0}} \frac{ \norm{Ce^{At} B^{-1} y}^2_2}{\norm{y}^2_2} \,,
\end{equation}
where the last term is the induced matrix norm $\norm{Ce^{At}}^2_2$ that can be computed by using SVD. 

An important issue must be considered, however, when we define these inner products associated with matrices $B$ and $C$. Let us assume that we measure only part of the system (e.g., the rotor speeds in the case where we are  concerned only with rotor stability). In this case the matrices $C$ and $B$ can be singular because the columns associated with the unmeasured state portion are zero. If $C$ is singular, however, applying the standard definition of the induced matrix norm is troublesome. A semidefinite matrix $C$ has an associated semidefinite inner product that induces a seminorm. Using the standard definition of induced matrix norm, we have
\begin{equation}
    \label{eq:seminorm_bad}
    \norm{A}_C = \maxx_{\substack{\norm{Cx}_2 \leq 1 \\ Cx \neq 0}} \norm{CAx}^2_2 \,.
\end{equation}
Given that $C$ is singular, every $x \neq 0$ that falls into the nullspace of $C$ can increase indefinitely and result in an unbounded value of the matrix norm. This issue has been considered in the non-normal dynamics literature, and the solution has been either to solve the problem via optimization \cite{Foures2012} or to define a normalized growth, where the nullspace of $C$ is replaced by the identity \cite{Blumenthal2016}. The former option is not desirable because the ability to employ the association between the largest singular value and the solution to the optimization problem is lost; the latter option is undesirable because it might result in solutions that are difficult to interpret. Instead, we propose a different approach that still relies on  SVD and does not require defining a normalized energy function. We make use of an alternative definition of the seminorm that, unlike \eqref{eq:seminorm_bad},  is not ill-posed. This definition, introduced in \cite{Jafarpour2022}, is
\begin{equation}
    \norm{A}_C = \maxx_{\norm{x}_C \leq 1 \,, \, x \in \mathcal{N}^\perp(C)} \norm{Ax}^2_C \,.
\end{equation}
Here the values of $x$ are restricted to lie in the subspace orthogonal to the nullspace of $C$, $\mathcal{N^\perp(C)}$, which prevents the norm from becoming unbounded. If we obtain a basis for $\mathcal{N(C)}$ and we arrange the vectors as the columns of a matrix $L$, we can enforce that all the vectors $x \in \mathcal{N}^\perp(C)$ by including the constraint $Lx = 0$. As in \eqref{eq:opt_problem}, the computation of the norm can be written as a constrained optimization problem:
\begin{equation}
\label{eq:opt_problem_constrained}
\begin{aligned}
& \underset{x(0)}{\text{maximize}}
& & \norm{x(t)}^2 \\
& \text{subject to}
& & \norm{x_0} = 1 \,,
& & Lx_0 = 0 \,.
\end{aligned}
\end{equation}
By focusing on this optimization problem, one can still  use  SVD to compute the seminorm. An equivalent optimization problem was solved by Golub using the QR decomposition of $L$ \cite{Golub1973}. In most of our cases, the matrix $C$ is constructed to weight a small subset of the state variables (e.g., the rotor speeds). In these cases, characterizing the nullspace and the subspace orthogonal to it is fairly trivial. We can obtain a basis for this $\mathcal{N}^\perp(C)$ and assemble the basis vectors as the matrix $F$. Then, we can make the change of variables $x = Fy$, which will ensure $x \in \mathcal{N}^\perp(C)$. Assume that our weighting matrix is concerned with only a subset of the state variables (equations). Then we can reorder the state variables (equations) such that the weighting matrix can be written
\begin{equation}
    \label{eq:semidef_weight}
    C = \begin{bmatrix} C_{11} & 0\\ 0 &  0 \end{bmatrix} \in \mathbf{R}^{n \times n} \,,
\end{equation}
where $C_{11} \in \mathbf{R}^{k \times k}$ is a nonsingular matrix. If we define a matrix
\begin{equation}
    F = \begin{bmatrix} I\\ 0 \end{bmatrix} \in \mathbf{R}^{n \times k} \,,
\end{equation}
where $I \in \mathbf{R}^{k \times k}$ is the identity matrix, then all nonzero vectors $x = Fz$ pertain to $\mathcal{N}^\perp(C)$. If we introduce another change of variables, $y = \hat{C}z$, where
\begin{equation}
    \hat{C} = \begin{bmatrix} C_{11}^{-1}\\ 0 \end{bmatrix} \in \mathbf{R}^{n \times k} \,,
\end{equation}
then
\begin{align}
    \norm{A}_C &= \maxx_{\norm{x}_C \leq 1 \,, \, x \in \mathcal{N}^\perp(C)} \norm{Ax}^2_C \,. \\
    &=  \maxx_{\norm{z}_2 = 1}  \frac{\norm{Ce^{At}L\hat{C}z}^2_2}{\norm{z}^2_2} \,.
\end{align}


This setup allows us to analyze the linearized system as a system with inputs and outputs. The linear dynamical system is
\begin{align}
    \begin{split}
        \dot{x} &= Ax + Bu \,, \\
        y       &= Cx \,,
    \end{split}
\end{align}
with $A \in \mathbf{R}^{n \times n}$, $B \in \mathbf{R}^{n \times m}$, and $C \in \mathbf{R}^{k \times n}$. Taking the Laplace transform of this system, we obtain the resolvent
\begin{equation}
    Y(s) = C(sI-A)^{-1}BU(s) \,.
\end{equation}

By considering $x(0) = 0$ and applying an impulse input we obtain the impulse response map
\begin{equation}
    \mathcal{H} = Ce^{tA}B \in \mathbf{R}^{k \times m} \,,
\end{equation}
which maps an input impulse at $t=0$ to the output at $t$. Then, to analyze the impulse perturbation that leads to the maximum growth in the output, we compute
\begin{equation}
    \maxx_{\substack{\norm{u(0)} \leq 1 \\ u(0) \neq 0}} \frac{ \norm{y(t)}^2_2}{\norm{u(0)}^2_2} = \maxx_{\substack{\norm{u(0)} \leq 1 \\ u(0) \neq 0}}  \frac{ \norm{Ce^{tA}B u(0)}^2_2}{\norm{u(0)}^2_2} =\norm{e^{At}}^2_{(C,B)} \,.
\end{equation}

Taking it all together, we can describe an algorithm to compute the maximum growth in a finite time span $t \in [0, T_{\textit{max}} ]$. This is shown in Algorithm \ref{alg:max_growth}. We can discretize this time interval in $k$ points, and for each point we compute the SVD of the matrix product $Ce^{At}B^{-1}$. We then seek the maximum singular value along this interval and the singular vector that indicates the perturbation direction with the maximum growth.

\begin{algorithm}
    \caption{MaxGrowth}
    \label{alg:max_growth}
    \begin{algorithmic}[1] 
        \Procedure{MaxGrowth}{}
        \State Given $A$, $B$, $C$, $dT$, $n$, initialize $s_{\textit{max}} = 1$
        \For{$k = 0, 1, 2, \cdots n $}
            \State $t = kdT$ 
            \State compute $Ce^{At}B^{-1} = U \Sigma V^\top $
            \State set $s_1 = \Sigma(1,1)$
            \State set $v_1 = V^\top (:,1)$
            \If{$s_1 > s_{\textit{max}}$}
            \State $s_{\textit{max}} = s_1$
            \State $v_{\textit{max}} = v_1$
            \EndIf{}
        \EndFor{}
        \State \Return $s_{\textit{max}}$, $v_{\textit{max}}$
        \EndProcedure
    \end{algorithmic}
\end{algorithm}

\subsection{Algorithm Implementation}
We have shown how the maximum singular value and vector of the map $e^{At}$ gives us the maximum preasymptotic growth magnitude and direction. For very large systems, however, forming this operator and computing the SVD can quickly overwhelm computer memory since the map is, in general, a dense matrix. This means that, for large systems, the algorithm is \textit{memory-bound} rather than \textit{compute-bound}. However, we do not need to compute all the singular values of the map nor form it explicitly. Indeed, one can think of computing the SVD as solving an equivalent eigenvalue problem:
\begin{equation*} 
    (e^{At})^\top  e^{At} v = \sigma^2 v \,.
\end{equation*}
The maximum eigenvalue (singular value) can be approximated by using an iterative solution such as the Lanczos iteration. These types of iterative algorithms do not require forming the operator explicitly, but only its \textit{action}. In our case, rather than explicitly computing the (possibly) dense matrix $e^{At}$, we only need to compute the matrix-vector operations $e^{At}v$ and $(e^{At})^*v$. We  observe that  the resulting matrix-vector operation $e^{At}v$ is equivalent to solving the initial value problem from $t=0$ to $t=T$:
\begin{equation}
    \dot{x} = Ax \,, \quad x(0) = v \,.
\end{equation}
In the same way, to obtain $(e^{At})^\top v$, we can solve the associated \textit{adjoint} system:
\begin{equation}
    \dot{y} = A^\top y \,, \quad y(0) = x(T) \,.
\end{equation}
These linear systems can be efficiently solved with explicit numerical methods that  require only sparse matrix-vector multiplications with the associated memory savings. For instance, the fourth-order Runge--Kutta method to integrate our linear dynamical system can be shown to reduce to
\begin{equation}
    x_{n + 1} = x_n + hAx_n + \frac{h^2}{2}A^2x_n + \frac{h^3}{6}A^3x_n + \frac{h^4}{24}A^4x_n \,,
\end{equation}
which is a truncated version of the matrix exponential series. 

We note that the matrix $A$ often arises from the linearization of a DAE that models the transient dynamics of the system. We write this DAE as
\begin{align*}
    \dot{x} &= f(x, y, t) \,,\\
    0 &= g(x, y, t) \,,
\end{align*}
where $x \in \mathbf{R}^n$ is the dynamic state vector and $y \in \mathbf{R}^m$ is the algebraic state vector. We then can linearize around an equilibrium point and reduce it to a linear dynamical system. We first linearize
\begin{align}
\begin{split}
    \label{eq:linear_dae}
    \Delta \dot{x} &= f_x \Delta x + f_y \Delta y \,,\\
    0 &= g_x \Delta x + g_y \Delta y \,,    
\end{split}
\end{align}
and, using the Schur complement, we reduce it to an equivalent linear dynamical system,
\begin{equation}
    \Delta \dot{x} = J_r \Delta x \,, \quad J_r = f_x - f_y(g_y^{-1})g_x \,.
\end{equation}
It is \textit{not} necessary to form the matrix $J_r$ explicitly. Here we can also adopt a matrix-free approach that will allow us to preserve sparsity. Because we  need only the \textit{action} of $J_r$, we can instead store $f_x$, $f_y$, $g_x$, $g_y^{-1}$ and compute a series of matrix-vector products to obtain $J_rx$ and $J^\top_rx$. This procedure is shown in Algorithms \ref{alg:mfree1} and \ref{alg:mfree2}.

\begin{algorithm}
    \caption{Matrix-free reduced Jacobian MATMULT}\label{alg:mfree1}
    \begin{algorithmic}[1] 
        \Procedure{MatMult}{$f_x$, $f_y$, $g_x$, $g_y$, $\Delta x$}
        \State u = $g_x \Delta x$
        \State y = $g_y^{-1} u$
        \State z = $f_y y$
        \State \Return $f_x \Delta x + z$
        \EndProcedure
    \end{algorithmic}
\end{algorithm}

\begin{algorithm}
    \caption{Matrix-free reduced Jacobian RMATMULT}\label{alg:mfree2}
    \begin{algorithmic}[1] 
        \Procedure{RMatMult}{$f_x$, $f_y$, $g_x$, $g_y$, $\Delta x$}
        \State u = $f^\top_y \Delta x$
        \State y = $(g^\top_y)^{-1} u$
        \State z = $g^\top_x y$
        \State \Return $f^\top_x \Delta x + z$
        \EndProcedure
    \end{algorithmic}
\end{algorithm}
Of course, the inverse $g_y^{-1}$ does not need to be formed explicitly. Instead, we can perform a sparse LU decomposition and store the factors to be reused each time.

%% file: sections/cases.tex
\subsection{Simplified Voltage Regulation Dynamics}
For the first example, we  analyze the linearized voltage regulation dynamics of a simplified power system. We use  a model introduced by Ilic and Zaborszky \cite[Chapter~7]{Ilic2000} that consists of a simplified differential-algebraic set of equations representing excitation dynamics. These equations are
\begin{align*}
\dot{e}' &= \frac{1}{T'_{d0}} \left( \frac{x_l + x_d}{x'}e' + \frac{x_d - x_d'}{x'}\frac{v^2 + x'Q}{e'} + e_{fd} \right) \,, \\
\dot{e}_{fd} &= \frac{1}{T} \left( -(e_{fd} - e_{fd}^0) - K \left(  \lambda \right) \right) \,, \\
0 &= e'^2 v^2 - (x' P)^2 - (x'Q + v^2)^2 \,, \\
\lambda &= \frac{\sqrt{(x_d')^2 + (xQ + v^2)^2}}{v} - e_r \,.
\end{align*}
In this model, the dynamic states are $e'$ and $v$, which are the generator voltage and the field voltage, respectively. The algebraic state $v$ is the voltage magnitude at the load terminal. For this experiment, we chose parameters
$T_{d0} = 10$, 
$T = 0.5$, 
$E_{fd0} = 1.6$, 
$E_r = 1.0$, 
$x_l = 0.1$, 
$x_d = 1.0$, 
$x_d' = 0.2$, 
$x' = 0.1$, 
$P = 0.5$, 
$Q = 0.2$, 
$H = 0.1$,  and
$B = 0.1$.

To study the effects of non-normality on transient preasymptotic growth, we examined two cases: (1) the gain of the excitation control $K$  is set to $0.5$ (System 1) and (2) the gain is set to $4.0$ (system 2). If we linearize the system and compute the implicit Jacobians, the results are
\begin{equation*}
    J_1 = \begin{bmatrix} -0.082 & 0.1\\ -1.015 & -2\\ \end{bmatrix}, \quad
    J_2 = \begin{bmatrix} -0.069 & 0.1\\ -8.123 & -2\\\end{bmatrix}\,.
\end{equation*}

\begin{figure}[h]
	\centering
	\includegraphics[width=0.45\textwidth]{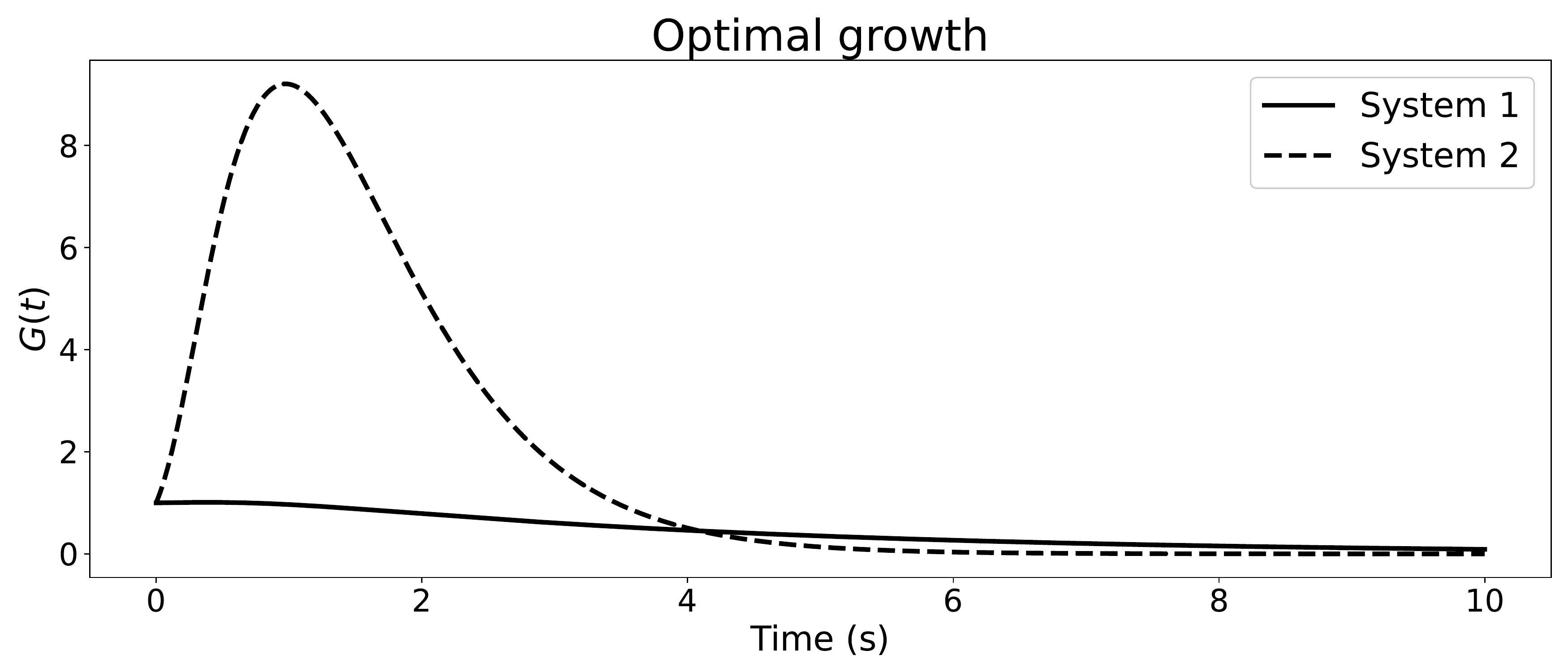}
	\caption{Optimal growth of systems 1 and 2 of case A computed with the algorithm}
	\label{fig:casea_optg}
\end{figure}

We observe that the bottom off-diagonal term is very sensitive to the changes in $K$. The spectra of these two matrices are $\sigma(J_1) = \{ -3.51, -1.94 \}$ and $\sigma(J_2) = \{ -4.24, -1.2 \}$, which in both cases indicate that the linearized system is asymptotically stable with eigenvalues of the same order of magnitude. However, the eigenvectors tell us another story. The eigenbasis for systems 1 and 2 are
\begin{equation*}
    V_1 = \begin{bmatrix} 0.878 & -0.053\\ -0.478 & 0.998\\ \end{bmatrix}, \quad
    V_2 = \begin{bmatrix} 0.159 & -0.076\\ -0.987 & 0.997\\ \end{bmatrix} \,.
\end{equation*}
We can see that in case 2 the angle between eigenvectors is smaller than in case 1. The condition number confirms this with $\kappa(V_1) = 1.79$ and $\kappa(V_2) = 23.82$. Using the algorithm introduced in Section \ref{section:methodology}, we can compute the optimal growth in both systems. This is plotted in Fig.~\ref{fig:casea_optg} and shows that system 2 exhibits a much more pronounced transient growth than does system 1 with a peak at $t=0.97$ seconds where $G(t)=9.2$.

Albeit a simple example, this experiment shows that eigenvalue magnitudes do not carry enough information to characterize preasymptotic growth.

\subsection{Classical Multimachine System}
\begin{figure}[h]
	\centering
	\includegraphics[width=0.45\textwidth]{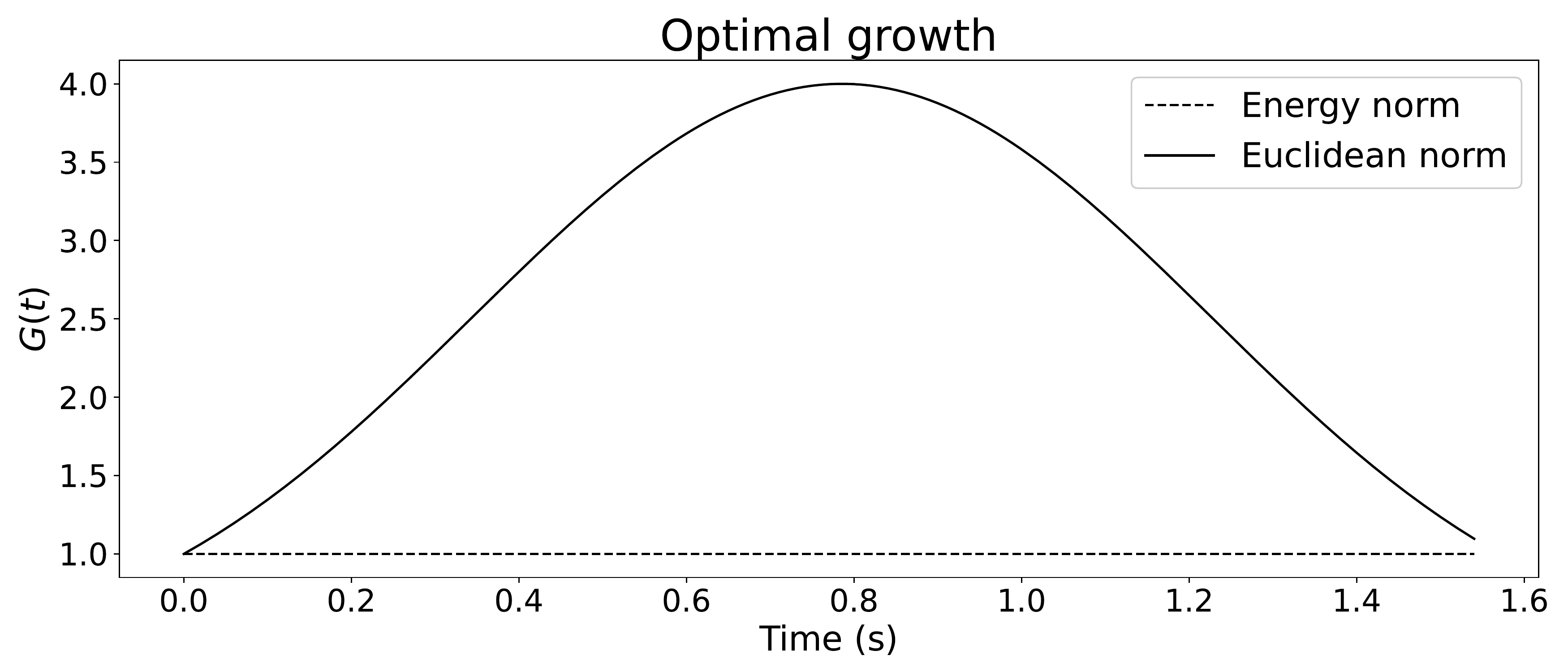}
	\caption{In continuous black: the optimal growth function using the Euclidian norm. This function is periodic.}
	\label{fig:caseb_growth1}
\end{figure}

\begin{figure}[h]
	\centering
	\includegraphics[width=0.45\textwidth]{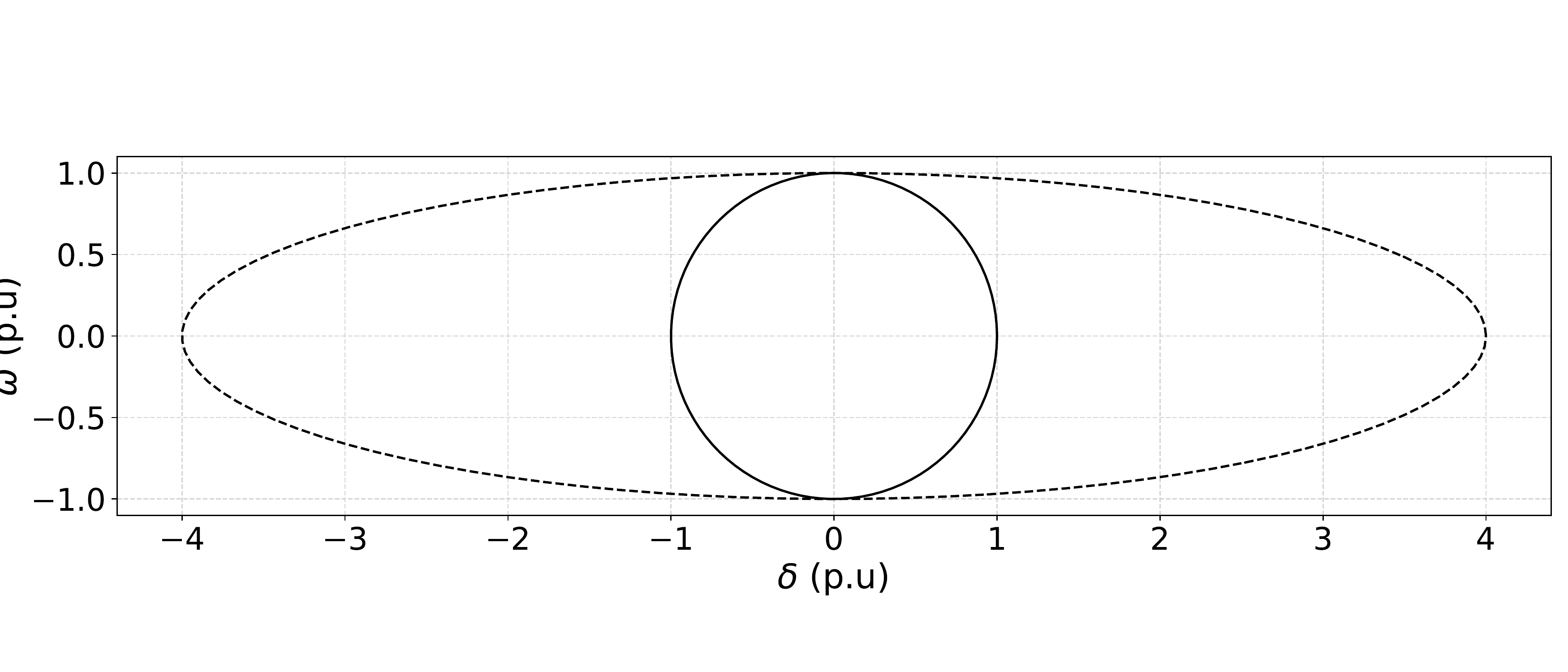}
	\caption{In continuous black: a unit ball of the state space such that $\norm{x}_2 \leq 1$. In discontinuous black:  a unit ball of the state space but using the energy norm instead, $\norm{x}_W \leq 1$}
	\label{fig:caseb_elp}
\end{figure}

The following example  illustrates the concepts of output norm. Recall that the classical model for power system dynamics can be expressed
 \begin{equation}
     H\ddot{\delta} + D\dot{\delta} = p_m - \Re(v\odot(Yv)^*) \,,
 \end{equation}
where $Y$ is the Kron-reduced admittance matrix. If we further simplify this system by assuming constant voltage magnitudes, the small angle hypothesis, and no transmission losses, we obtain a second-order linear differential equation:
\begin{equation}
    H\ddot{\delta} + D\dot{\delta}  + K\delta = p_m  \,.
\end{equation}
Writing this equation as a first-order system results in the system in Eq.~\eqref{eq:classical_linear}. As an example, we consider a two-machine, two-bus system, with no damping and no forcing term. We set the inertia and reactance matrix to
\begin{equation}
    H = \begin{bmatrix} 1 & \\ & 1 \end{bmatrix}, \quad K = \begin{bmatrix} 2 & -2   \\ -2 & 2 \end{bmatrix} \,.
\end{equation}
Because we do not have a reference angle, the system has an extra degree of freedom that makes the matrix $K$ positive semi-definite. To avoid this situation (which would complicate our analysis), we can take the first machine as reference by introducing an additional constraint.  The two-machine system then can be described as an equivalent one-machine system. The operator in \eqref{eq:classical_linear} becomes
\begin{equation}
A = \begin{bmatrix}
  0 & 1\\
  -4 & 0
\end{bmatrix} \,.
\end{equation}
The matrix operator $A$ is not normal since  $A^TA \neq I$, which would indicate the possibility of non-normal growth. However, we also know that because of the lack of damping, this system is conservative and will not dissipate energy. We  compute the optimal energy growth using both the Euclidean and energy norm in \eqref{eq:energy_norm}. The results of this computation are shown in Fig. \ref{fig:caseb_growth1}, where we can see that the Euclidean norm leads to a growth that peaks around $t=0.8$ seconds. Using the energy norm, however, we observe that the growth function is constant at 1, which indicates no energy change.
Recall that, for a conservative oscillator, the trajectories remain in an ellipsoid of constant energy, as described by
\begin{equation}
\label{eq:energy_ellipsoid}
    \mathcal{E}(x) = \mathcal{E}(\delta, \omega) = \frac{1}{2}\delta^T K \delta +  \frac{1}{2}\omega^T H \omega \,.
\end{equation}
This means that for an initial condition $x_0$, accounting for energy conservation, $\mathcal{E}(x_0) = \mathcal{E}(e^{At}x_0)$. This can also be written by using the energy norm $\norm{x_0}_W = \norm{e^{At}x_0}_W$. Hence we can write
\begin{equation}
   \maxx_{\substack{\norm{x}_W \leq 1 \\ x \neq 0}}  \frac{ \norm{e^{At} x}^2_W}{\norm{x}^2_W} = \maxx_{\substack{\norm{x}_W \leq 1 \\ x \neq 0}}  \frac{ \norm{x}^2_W}{\norm{x}^2_W} = 1 \,.  
\end{equation}
The energy is preserved. When we use the Euclidian norm, however, we observe transient energy growth that peaks at $t=0.8$ s. In Fig. \ref{fig:caseb_elp}, we draw a diagram to help understand this phenomenon. The growth maximization problem (from the Euclidean perspective) consists of finding a point in the unit ball such that the energy at time $t$ is maximized. Because of \eqref{eq:energy_ellipsoid}, we know that the system describes constant energy ellipsoids. For the given unit ball of Fig. \ref{fig:caseb_elp}, the maximum energy ellipsoid is obtained when $\omega_0 = \pm 1$. After some time, the kinetic energy is transformed into potential energy such that $\delta = \pm 4$. This transfer from kinetic to potential energy explains the numerical growth observed with the Euclidean norm.

With this experiment, we can see the importance of choosing an adequate output norm. In this particular case, the Euclidean norm confers values of $\delta$ and $\omega$ the same weight whereas the energy norm weights them in such a way that transforms the ellipsoid into a circle. This does not mean that using the Euclidean norm results in incorrect results but, instead, that the interpretation of the results depends on the measure.  This reliance on a measure is a distinctive trait our the non-modal technique.

\subsection{Two-Machine System with Classical Generator Model}
An additional level of detail is provided by the classical electromechanical model, which represents generator dynamics as a voltage source behind an impedance. In this section we consider an instance of this model introduced in \cite{Bretas2003}. The system consists of two generators connected to an infinite bus, and it can be described in ordinary differential equation form as
\begin{align*}
\dot{\delta}_1 &= \omega_1 \,, \\
\dot{\delta}_2 &= \omega_2 \,, \\
H_1\dot{\omega}_1 &= P_1 - \bar{B}_{10} \sin \delta_1 - \bar{G}_{10} \cos \delta_1 - \bar{B}_{12} \sin(\delta_1 - \delta_2) \\
&- \bar{G}_{12} \cos(\delta_1 - \delta_2) - D_1 \omega_1 \,, \\
H_2\dot{\omega}_1 &= P_2 - \bar{B}_{20} \sin \delta_2 - \bar{G}_{20} \cos \delta_2 - \bar{B}_{12} \sin(\delta_2 - \delta_1) \\
&- \bar{G}_{12} \cos(\delta_2 - \delta_1) - D_2 \omega_2 \,,
\end{align*}
where bus $0$ is the infinite bus, $\bar{B}_{ij} = v_iv_jB_{ij}$, $\bar{G}_{ij} = v_iv_jG_{ij}$, and $P_i = P_{\textit{mech}} - v_i^2 G_i$. Linearizing about a steady-state point, the Jacobian matrix has the same structure as \eqref{eq:classical_linear} but with the matrix $K$ being the partial derivatives of the third and fourth equations with respect to $\delta_1$ and $\delta_2$. If we set $D_1 = D_2 = 0$,  denote $\Omega = H^{-1}K$, and  take $W = \begin{bmatrix} 0 & I\end{bmatrix} \in \mathbf{R}^{2 \times 4}$, then we see that
\begin{align*}
    W e^{At} W^T &= I - \frac{1}{2!}\Omega t^2 + \frac{1}{4!} \Omega^2 t^4 - \frac{1}{6!} \Omega^3 t^6 + \dots \\
    &= \cos ( \sqrt{\Omega} t) \,.
\end{align*}
Here $\sqrt{\cdot}$ and $\cos(\cdot)$ are the \textit{matrix} square root and cosine.  We can also compute the energy growth as 
\begin{equation}
   \maxx_{\substack{\norm{x(0)}_2 \leq 1 \\ x(0) \neq 0}}  \mathcal{G}(t) = \maxx_{\substack{\norm{x(0)}_2 \leq 1 \\ x(0) \neq 0}} \frac{\norm{\cos(\sqrt{\Omega} t)x(0)}}{\norm{x(0)}} \,.
\end{equation}

\begin{figure}[h]
	\centering
	\includegraphics[width=0.45\textwidth]{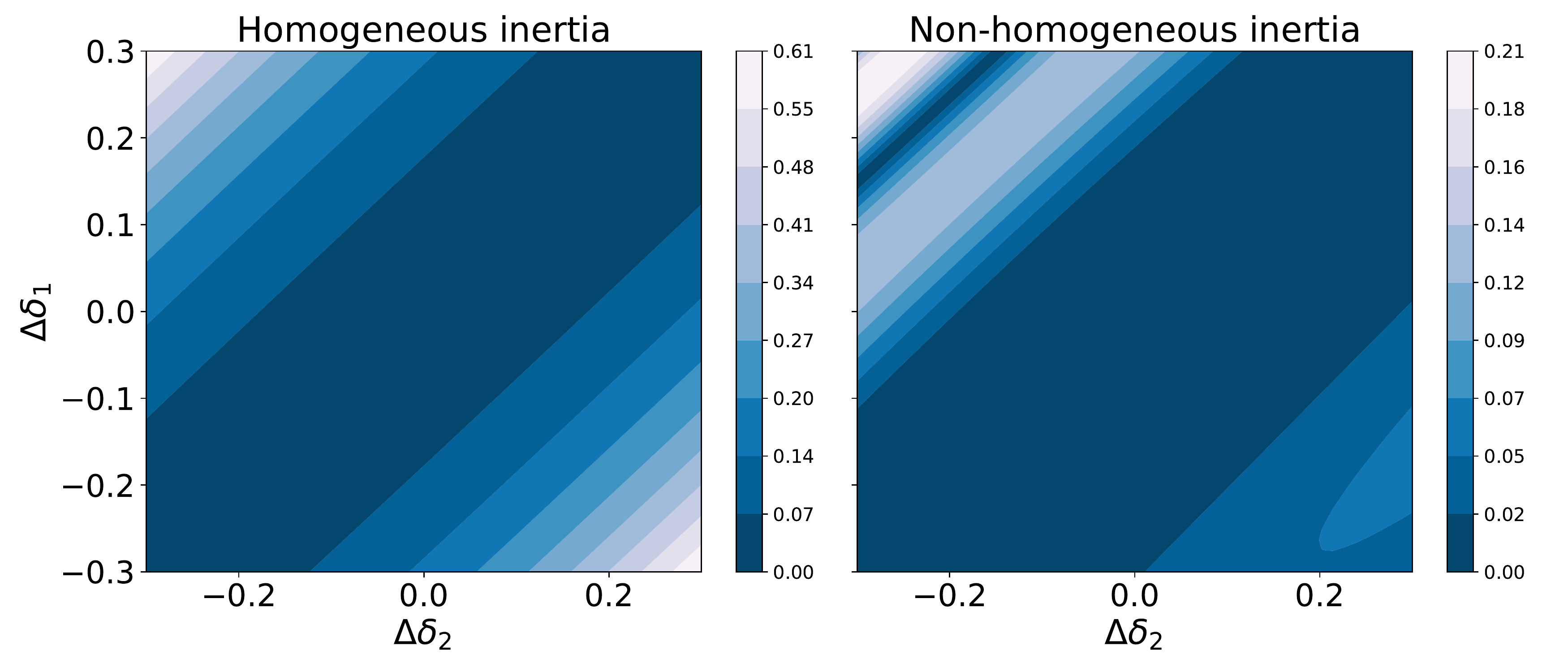}
	\caption{Henrici number for different angle deviations. On the left is a system with homogeneous inertia. We observe how the Henrici number increases with the increase of angle difference. On the right is a system with nonhomogeneous inertia where the dependence of the Henrici number with the angle difference is not as clear.}
	\label{fig:henrici}
\end{figure}

Notice that if we use the subordinate norm, $\norm{\cdot}_H$, then
\begin{align*}
    \norm{\cos(\sqrt{\Omega} t)}_H &= \norm{H^{\frac{1}{2}} \cos(\sqrt{\Omega} t)  (H^{\frac{1}{2}})^{-1}}_2 \\
    &= \norm{\cos(\sqrt {(H^{\frac{1}{2}})^{-T} K (H^{\frac{1}{2}})^{-1}}  t) }_2 \,,
\end{align*}
where $(H^{\frac{1}{2}})^{-T}$ is the inverse transpose of the Cholesky factor of $H$. If $K$ is symmetric, then the matrix product inside the square root is symmetric, which by the Cayley--Hamilton theorem indicates that we can find an orthogonal eigenbasis and the preasymptotic growth is determined completely by the eigenvalues, which are all less than 1 in magnitude---there is no increase of kinetic energy. When the losses of the network are negligible, the matrix $K$ takes the form
\begin{equation*}
    \begin{bmatrix}
    \bar{B}_{10} \cos{\delta_1} + \bar{B}_{12} \cos{(\Delta \delta_{12})} & -\bar{B}_{12} \cos{(\Delta \delta_{12})} \\
    -\bar{B}_{12} \cos{(\Delta \delta_{21})} & B_{20} \cos{\delta_2} + \bar{B}_{12} \cos{(\Delta \delta_{21})}
    \end{bmatrix} \,.
\end{equation*}
Because the cosine is even, the matrix is symmetric. This property can be broken when considering network losses that can lead to non-normal preasymptotic growth. In this case, the degree of non-normality can increase with the difference between $\delta_1$ and $\delta_2$. The Henrici index quantifies the departure from non-normality of a matrix $A$ and is defined as
\begin{equation}
    \nu(A) = \left( \norm{A}_F^2 - \sum_{j=1}^n |\lambda_j|^2 \right)^{\frac{1}{2}} \,.
\end{equation}
With this number and $\Omega$, we can compute the deviation of normality as the voltage angle of bus 1 and bus 2 varies. In Fig.~\ref{fig:henrici} we show this experiment in a case in which $H_1 = H_2$ and another in which $H_1 = 10H_2$. In the first one, we see how the increase in non-normality is directly proportional to the increase in voltage angle difference. In the case in which the inertia is nonhomogeneous, this relationship becomes more complicated.

Thus, in this example, we can see that non-normality can arise because of the network topology. Furthermore, for the classical electromechanical model, we can say that a lossy network is a prerequisite for the emergence of non-normality (independent of inertia).

\subsection{IEEE 9-Bus and 39-bus Systems}
In this subsection we examine the sub-asymptotic growth in the IEEE 9-bus and 39-bus systems using a positive sequence model with detailed synchronous generators and control. In our example we simulate the case using the Sauer--Pai model for the synchronous generators as well as the IEEESGO model for governor control \cite{Sauer1997}. The swing equation for the synchronous generator model can be written
\begin{equation}
    \dot{\omega} = \frac{1}{H}\left(\frac{p_m - D\omega}{1 + \omega} - \psi_{de}i_q + \psi_{qe}i_d\right),
\end{equation}
where $i_d$, $i_q$ are the stator currents and $\psi_{de}i_q$, $\psi_{de}i_d$ are the stator fluxes. When a governor is present, $p_m$ will depend on the rotor frequency. Because we now use  more realistic synchronous machine models with control, the structural simplicity of the models in the preceding sections is lost. One can still, of course, apply the framework introduced in this paper, but it is harder to determine the normality or non-normality of the operator of interest. 

The resulting DAE system is linearized according to \eqref{eq:linear_dae}. We can construct the matrix $C$ such as in \eqref{eq:semidef_weight}, where $C_{11}$ is a square diagonal matrix with a dimension corresponding to the number of machines and containing the square root inertia of each machine in the diagonal entries. We can then compute the optimal growth of the weighted machine frequency deviations.

In our experiment we take the base scenario and increase the loading by a factor, $\alpha$, of $1.6$ and $2.2$. If we compute the reduced Jacobian of the system and examine the three eigenvalues with the minimum absolute real part, we see that these are of similar magnitude for the three loading conditions (Table \ref{table:ieee9}). As we already know, these eigenvalues do not characterize the dynamics of the system completely, and the normality of the linear operator plays an important role. 
The Henrici index is computed for these three loading cases as shown in Table \ref{table:ieee9} and clearly shows that as we increase the loading, the non-normality of the operator increases. 

\begin{table}
	\renewcommand{\arraystretch}{1.3}
	\caption{Spectra closest to the origin and Henrici's departure from normality index. Case IEEE9}
	\label{table:ieee9}
	\centering
	\review{
	\begin{tabular}{c c c c c}
		\hline
		\bfseries $\alpha$ & $\lambda_1$ &  $\lambda_2 $ & $\lambda_3 $ & $\nu$  \\
		\hline
		1.0 & -4.99E-02 & -4.99E-02 & -5.85E-02 & 0 \\
		1.6 & -4.99E-02 & -4.99E-02 & -5.87E-02 & 0.277 \\
		2.2 & -4.26E-02 & -4.99E-02 & -5.02E-02 & 1.335 \\
		\hline
	\end{tabular}}
\end{table}

\begin{figure}[ht]
	\centering
	\includegraphics[width=0.45\textwidth]{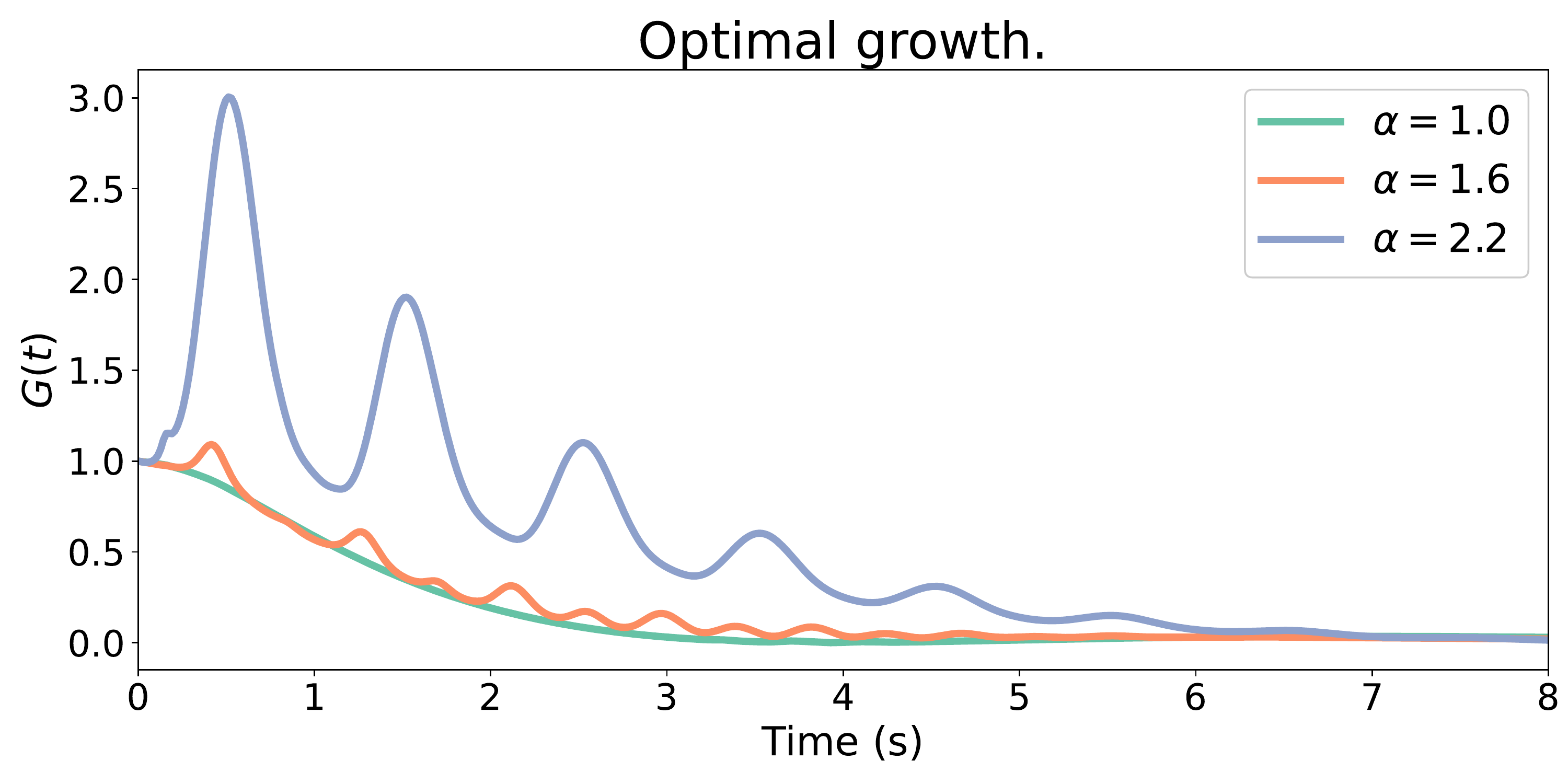}
	\caption{Optimal growth of the 9-bus system for increasing loading conditions}
	\label{fig:ieee9_1}
\end{figure}

\begin{figure}[ht]
	\centering
	\includegraphics[width=0.45\textwidth]{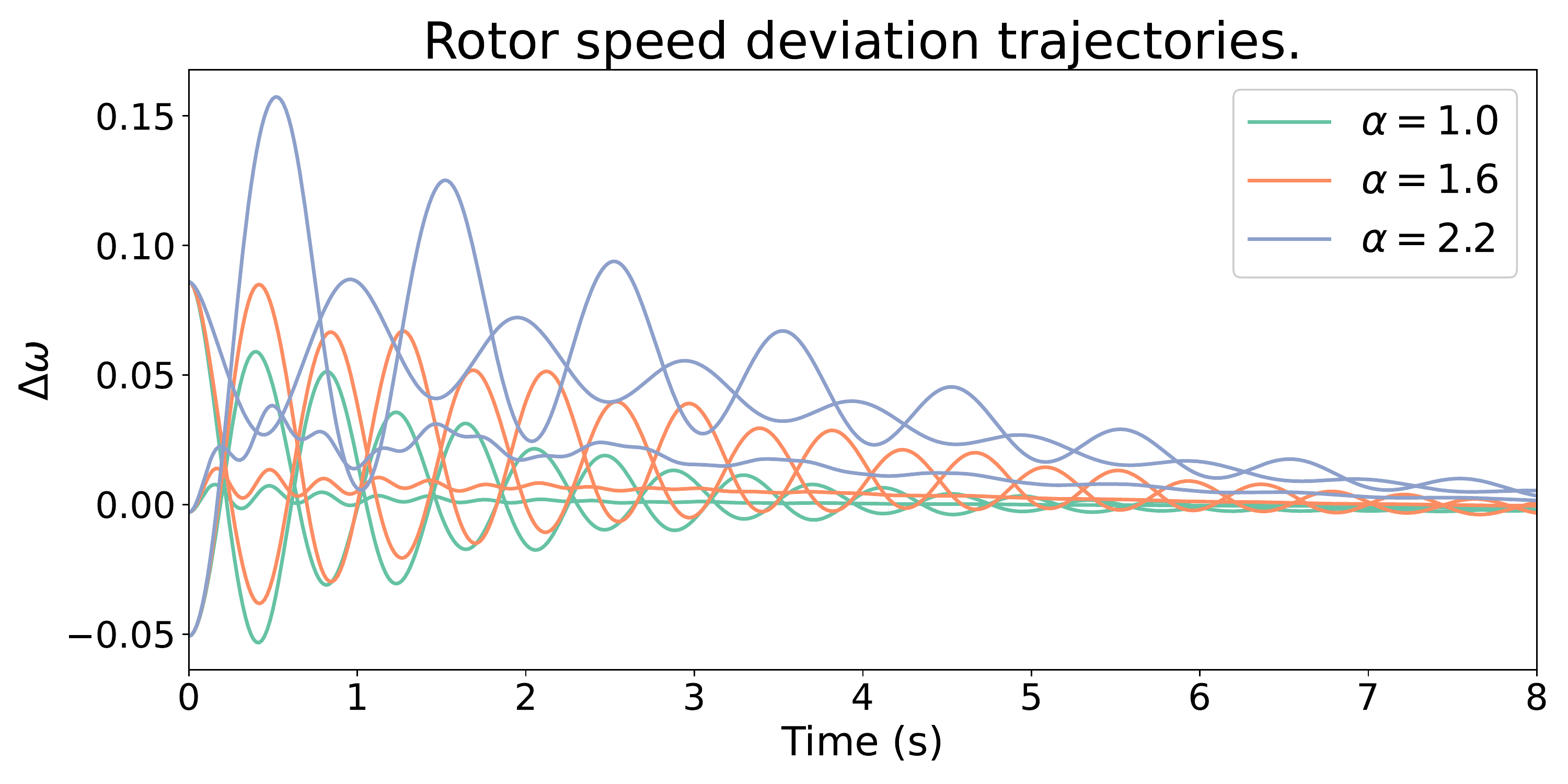}
	\caption{Impulse response of the 9-bus system for maximum growth perturbation}
	\label{fig:ieee9_2}
\end{figure}

We then compute the maximum growth using Algorithm \ref{alg:max_growth}. The results are shown in Fig.~\ref{fig:ieee9_1}. For the base case ($\alpha = 1$), the maximum growth peaks at the beginning with a relative growth of 1. This means that the energy is dissipated. As we increase the loading, the system begins to display a relatively large subtransient growth before dissipating the energy. In the case of $\alpha=2.2$, this occurs at $t=0.513$ sec.; and, from the SVD analysis, we determine that the perturbation that leads to optimal growth is
\begin{equation*}
 v_1 = \begin{bmatrix}
 0.86 & -0.03 & -0.508 
 \end{bmatrix}  \,.
\end{equation*}
In Fig.~\ref{fig:ieee9_2} we plot the response of the system with increasing loading conditions to the optimal perturbation vector.  We can see that the sensitivity of the $\alpha=2.2$ system to the same perturbation is much higher. The speed deviation of generator 3 overshoots and peaks at $t = 0.513$, which aligns with the growth analysis. Furthermore, we can see that, whereas the response of the $\alpha=1.6$ system seems a \textit{perturbed} version of the $\alpha=1.0$ system (i.e., the trajectories differ numerically but exhibit similar behavior), the response of system $\alpha = 2.2$ seems qualitatively different.

\begin{figure}[ht]
	\centering
	\includegraphics[width=0.45\textwidth]{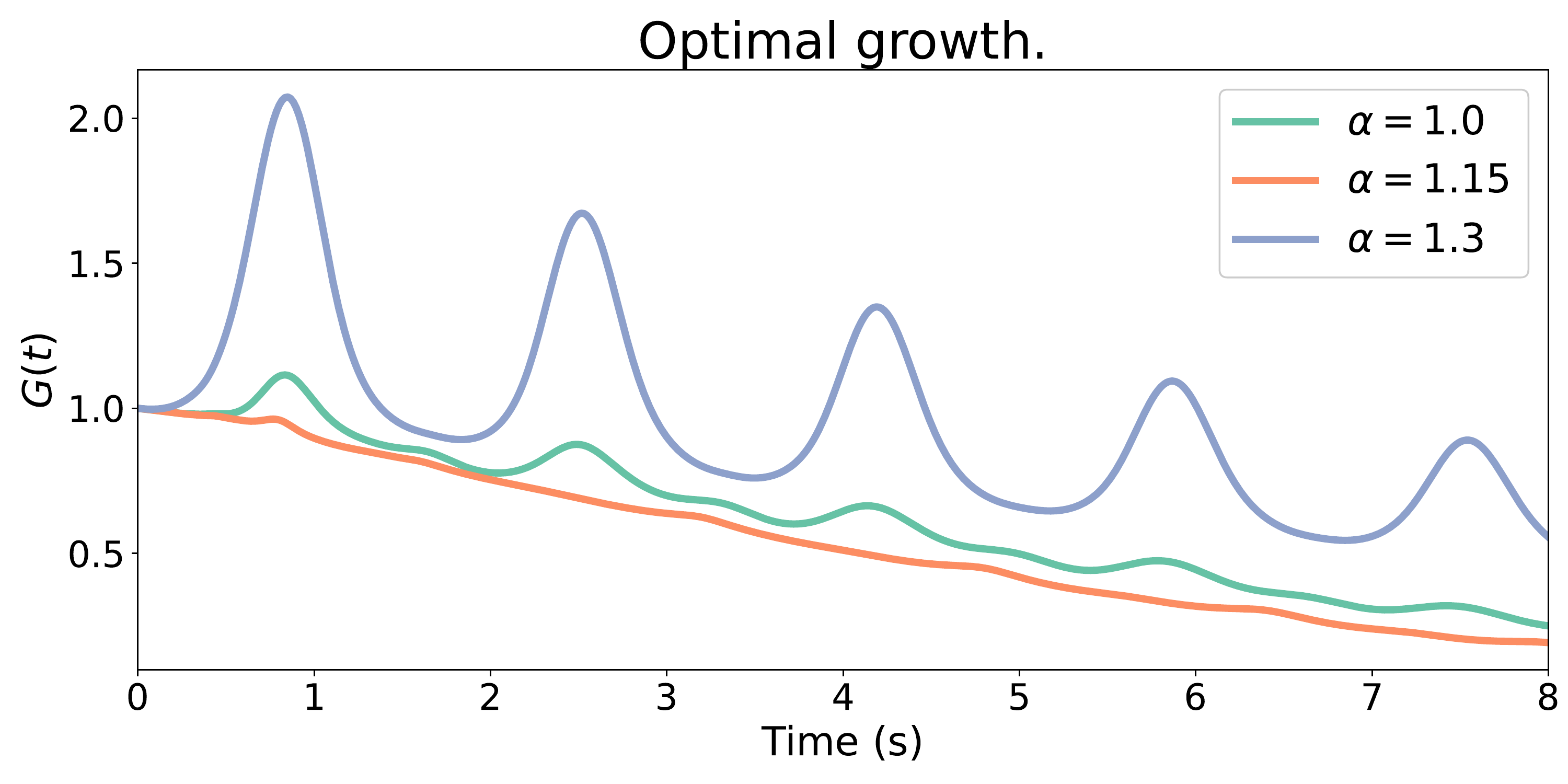}
	\caption{Optimal growth of the 39-bus system for increasing loading conditions}
	\label{fig:ieee_ne}
\end{figure}

\begin{figure}[h]
	\centering
	\includegraphics[width=0.45\textwidth]{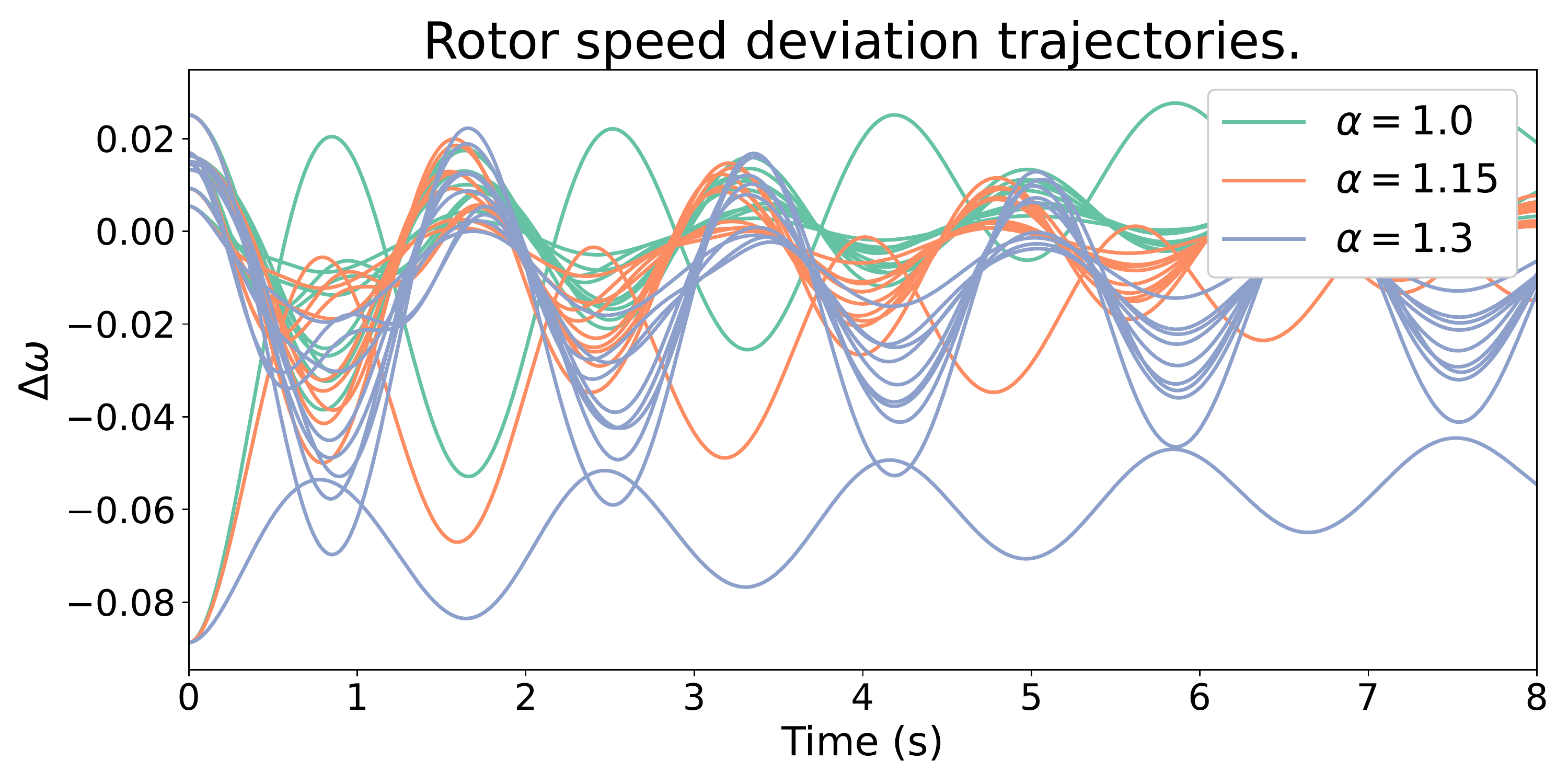}
	\caption{Impulse response of the 39-bus system for maximum growth perturbation}
	\label{fig:ieee_ne2}
\end{figure}

We now perform the same analysis on the 39-bus system. We increase the loading conditions up to a factor of $1.3$, and we observe its effects on the spectra and the normality. In Table \ref{table:caseNE} we see that, as in the previous system, the spectra of the Jacobian matrix are relatively similar. In the same way as for the 9-bus system, however, the Henrici number reveals that these systems will behave differently. This is corroborated in Fig.~\ref{fig:ieee_ne} where we see that the system with $\alpha=1.3$ exhibits transient growth peaking around $t = 0.9$. The full response to the optimal perturbation is shown in Fig.~\ref{fig:ieee_ne2}, where we can see how, for $\alpha=1.3$, the speed deviation of one of the generators has difficulty reaching equilibrium. An additional insight, over the 9-bus system, is that non-normality does not necessarily grow with the increase of loading conditions. As we can see in Table \ref{table:caseNE}, the Henrici number for $\alpha=1.0$ is higher than for $\alpha=1.15$. This points to the fact that, as we saw in the preceding subsection, non-normality is associated with the angle difference but this relationship is not necessarily linear.   

\begin{table}
	\renewcommand{\arraystretch}{1.3}
	\caption{Spectra closest to the origin and Henrici's departure from normality index in the 39-bus system}
	\label{table:caseNE}
	\centering
	\review{
	\begin{tabular}{c c c c c}
		\hline
		\bfseries $\alpha$ & $\lambda_1$ &  $\lambda_2 $ & $\lambda_3 $ & $\nu$  \\
		\hline
		1.0 & -4.99E-02 & -4.99E-02 & -4.99E-02 & 0.295 \\
		1.15 & -4.99E-02 & -4.99E-02 & -4.99E-02 & 0.0 \\
		1.3 & -3.99E-02 & -4.74E-02 & -4.98E-02 & 0.836 \\
		\hline
	\end{tabular}}
\end{table}

\subsection{Scalability to Large Systems}
The scalability of the matrix-free method is also tested for the large-scale systems to compute the exponential matrix operator norm described in Section \ref{section:methodology}. The experiment is performed on a workstation equipped with an Intel Xeon E5-2650 v2 processor and 16GB of DDR2 DRAM. We use the Python programming language and, in particular, sparse linear algebra libraries such as the UMFPACK wrapper. We implement Algorithm \ref{alg:max_growth} both with an explicit computation of the exponential map and with a matrix-free approach. We compute 5 points in time in an interval such that $t \in (0.0, 1.0) $. We perform this experiment for 7 test systems of an increasing number of buses (9, 39, 200, 2000, 10k, 20k, 70k). The results are shown in Table~\ref{table:scalability}. We can see that while computing the explicit map is faster, for larger cases the memory requirements start to grow quickly and the computation becomes troublesome for very large systems. Indeed, for the ACTIVSg10k case, the memory requirements of the explicit method are two orders of magnitude above the matrix-free method. For ACTIVSg25k and larger, the program runs out of memory. In this case we can still approximate the memory consumption using the dimension of the linear map and the size of a double-precision floating-point number in memory. The storage requirements increase quadratically with the size of the matrix; and, as a result, the scalability of the explicit method suffers for larger systems. In contrast, the matrix-free method is slower than the explicit method, but it can scale. Thus, we can observe how using the matrix-free method is a requisite for using the algorithms to study  very large systems.

\begin{table}
	\renewcommand{\arraystretch}{1.3}
	\caption{Scalability of the explicit and matrix-free algorithms. In parentheses are the matrix-free results.}
	\label{table:scalability}
	\centering
	\review{
	\begin{tabular}{c r c c}
		\hline
		\bfseries System & \# of buses & Time (s) &  Memory (MB) \\
		\hline
		IEEE9 &  9 & 4.8E-02 (1.7E-01) & 4.4E+01 (8.6E+00) \\
		New England & 39 & 6.1E-02 (5.9E-01) & 2.8E+00 (8.6E+00)  \\
		ACTIVSg200 & 200 & 6.5E-01 (3.4E+00) & 1.1E+01 (2.5E+00)\\
		ACTIVSg2000 & 2,000 & 3.1E-00 (3.8E+00) & 4.4E+01 (8.6E+00)  \\
		ACTIVSg10k & 10,000 & 2.8E+02 (1.1E+02) & 1.0E+02 (2.0E+00) \\
		ACTIVSg25k & 25,000 & FAIL (2.8E+02) & FAIL (1.8E+01) \\
		ACTIVSg70k & 70,000 & FAIL (8.8E+02) &  FAIL (3.5E+01) \\
		\hline
	\end{tabular}}
\end{table}